\documentclass[a4paper,10pt]{article}
\usepackage[utf8x]{inputenc}
\usepackage{amssymb}

\title{Constructing  Mironov cycles in complex Grassmanians }
\author{ Nikolai Tyurin\footnote{The author is partially supported by Laboratory of Mirror Symmetry NRU HSE, RF Government grant, ag. N 14.641.31.0001}\\
BLTPh JINR (Dubna) and NRU HSE (Moscow)}

\date{}

\begin{document}

\maketitle

A. Mironov in paper [1] proposed a construction of lagrangian submanifolds in $\mathbb{C}^n$ and $\mathbb{C} \mathbb{P}^n$; there he was mostly motivated by the fact that
these lagrangian submanifolds (which can have in general self intersections, therefore below we call them lagrangian cycles) present new example of minimal or Hamiltonian minimal lagrangian submanifolds. However the Mironov construction of lagrangian cycles  itself can be directly extended to much  wider class of compact algrebraic varieties: namely it works in the case when algebraic variety $X$ of complex dimension $n$ admits $T^k$ - action and an anti - holomorphic involution such that the real part $X_{\mathbb{R}} \subset X$ has real dimension $n$ and is transversal to the torus action. For this case, as we show in [2], one has families of lagrangian submanifolds and cycles.

In the present small text we show how the construction of Mironov cycles works for the complex Grassmannians, resulting in simple examples of smooth lagrangian submanifolds
in ${\rm Gr}(k, n+1)$, equipped with a standard Kahler form under the Pl\"{u}cker embedding. For sure the text is not complete but in the new reality we would like to fix it, hoping to continue the investigations and to present in a future complete list of Mironov cycles in ${\rm Gr}(k, n+1)$.

{\bf Acknowledgements.} The author cordially thanks A. Kuznetsov and P. Pushkar' for valuable discussions and remarks.

$$$$

{\bf General theory.} Let $(X, \omega)$ be a simply connected compact smooth  symplectic manifold of real dimension $2n$. Suppose that it admits an incomplete toric action so
there are moment maps $f_1, ..., f_k$ which commute with respect to the standard Poisson brackets and generate Hamiltonian $T^k$ - torus action on $M$. Take a generic set of values
$c_1, ..., c_k$ such that the common level set $N(c_1, ..., c_k) = \{ f_i = c_i, i = 1, ..., k\}$ does not intersect the determinantal locus $\Delta(f_1, ..., f_k) =
\{ X_{f_1} \wedge ... \wedge X_{f_k} = 0 \}$ and contains an isotropical $n-k$ - dimensional submanifold $S_0 \subset N(c_1, ..., c_k)$ which is transversal at each point to the torus action which means that $T_p S_0$ and $<X_{f_1}(p), ..., X_{f_k}(p)>$ are transversal in $T_p N(c_1, ..., c_k)$ at each point $p \in S_0$.

Then the toric action applied to $S_0$ generates a real $n$ - dimensional cycle $T^k(S_0)$ which is lagrangian; depending on the situation it can be smooth or
have self intersections, moreover in some cases one can take even special values of $f_i$ and nevertheless get smooth lagrangian submanifolds. The proof is based on the fact that
in the present situation the linear span $<T_p(c_1, ..., c_k), X_{f_1}(p), ..., X_{f_k}(p)>$ is lagrangian subspace in $T_p M$ for each $p \in S_0$, and the lagrangian condition is stable under the  Hamiltonian  $T^k$ - action (we present the details and  the simplest examples in [2]).
 
{\bf Algebraic varieties.} Any compact algebraic variety $X$ can be considered as a real symplectic variety: by the very definition (see f.e. [3]) any $X$ admits a very ample line bundle
$L \to X$ such that the corresponding complete linear system $\vert L \vert$ generates an embedding $\phi_L : X \hookrightarrow \mathbb{C} \mathbb{P}^N$, and the lifting
$\omega_L = \phi_L^* \Omega_{FS}$ of a standard Kahler form $\Omega_{FS}$ of the Fubini - Study metric gives a symplectic (Kahler) form $\omega_L$ on $X$. This form is not unique,
but the cohomology class $[\omega_L] = c_1(L) \in H^2(X, \mathbb{Z})$ if fixed, and one can expect that the lagrangian geometry of $(X, \omega_L)$ depends on this cohomology class only.
Note however that for different very ample line bundles over a fixed $X$ one can expect rather different lagrangian geometries. 

Suppose that $(X, \omega_L)$ admits a Hamiltonian toric action; for our aims we formulate it as follows. Fix a homogenous coordinate system $[z_0: ... : z_N]$ on $\mathbb{C} \mathbb{P}^N$, compatible with fixed Fubini - Study form $\Omega_{FS}$, and consider the standard moment maps of the form
$$
F_i = \frac{\sum_{j=0}^N \lambda_{ij} \vert z_j \vert^2}{\sum_{j = 0}^N \vert z_j \vert^2}, \quad i = 1, ..., N, \quad \lambda_{ij} \in \mathbb{Z},
$$
such that if we add row $(\lambda_{0j}) = (1, ..., 1)$ then integer valued matrix $\Lambda = ( \lambda_{ij} ), 0 \leq i, j  \leq N,$ is non degenerated.
Then we suppose that it exists such $\Lambda$ that $k$ moment maps $F_i$ (without loss of generality we can think that they are $F_1, ..., F_k$) preserve by the Hamiltonian action
the image $\phi_L(X) \subset \mathbb{C} \mathbb{P}^N$, and these $F_i$'s generate the corresponding $T^k$ - action.

For certain algebraic varieties it is possible to find isotropical submanifolds in the common level sets $N(c_1, ..., c_k)$ moreless automatically. Suppose additionally that our $X$ admits
an appropriate anti holomorphic involution. Again we reformulate it in our simple terms: for the fixed coordinate system $[z_0: ...: z_N]$ consider
the map
$$
\sigma: [z_0: ... : z_N] \mapsto [\bar z_0: ... : \bar z_N];
$$
suppose that $\phi_L(X)$ is {\it real} with respect to $\sigma$ which means that $\sigma(\phi_L(X)) = \phi_L(X)$ and moreover that the real part
$X_{\mathbb{R}} = \{ x \in X \quad | \quad  \phi^*_L \sigma (x) = x \} \subset X$ is smooth real submanifold of dimension $n$.

Then $X_{\mathbb{R}} \subset X$ is lagrangian with respect to $\omega_L$; consequently the intersection $ S_{\mathbb{R}}(c_1, ..., c_k) = X_{\mathbb{R}} \cap N(c_1, ..., c_k)$ is isotropical, and moreover if the value set $(c_1, ..., c_k)$ is generic then components of $S_{\mathbb{R}}(c_1, ..., c_k)$ are smooth and transversal to the $T^k$ - action. Therefore
we can apply General theory which leads to the construction of lagrangian cycles $T^n(S_{\mathbb{R}}(c_1, ..., c_k)$ in $(X, \omega_L)$. 

Indeed, the transversality of the toric action and the real part can be directly checked in $\mathbb{C} \mathbb{P}^N$:
the image of the real part $\phi_L(X_{\mathbb{R}}) = \phi_L(X) \cap \mathbb{R} \mathbb{P}^N$, and the transversality is clear from the coordinate description
of $T^k$ - action and the real part $\mathbb{R} \mathbb{P}^N \subset \mathbb{C} \mathbb{P}^N$.

Since this construction is  a natural extension of the construction for $\mathbb{C}^n$ and $\mathbb{C} \mathbb{P}^n$, presented in [1],
we call the resulting  $T^k (S_{\mathbb{R}}(c_1, ..., c_k)) \subset (X, \omega_L)$  {\it Mironov cycles} (or Mironov submanifolds
in the case when they are smooth).

{\bf Grassmannians.} The situation, presented above, takes place for  complex Grassmanians ${\rm Gr}(k, n+1)$ under the Pl\"{u}cker embedding to $ \mathbb{P}(\wedge^k \mathbb{C}^{n+1})$
(details on the geometry of ${\rm Gr}(k, n+1)$ can be found in [3]). In the discussion below we would like to avoid the algebraic machinery and construct certain Mironov cycles in ${\rm Gr} (k, n+1)$ using pure geometrical arguments, however all steps  can be explicitly  checked in the Pl\"{u}cker coordinates.

First or all we fix a hermitian structure on the source $\mathbb{C}^{n+1}$, compatible with coordinate system $(Z_0, ..., Z_n)$; this gives the corresponding K\"{a}hler structure on
the projective space $\mathbb{C} \mathbb{P}^n$, equipped with homogenous coordinates $[z_0: ... : z_n]$. Then the standard toric action generated by moment maps
$$
\mu_i = \frac{\vert z_i \vert^2}{\sum_{j=0}^n \vert z_j \vert^2}, \quad i = 1, ..., n,
$$
can be naturally extended to the spaces of all projective subspaces of $\mathbb{C} \mathbb{P}^n$ since evidently every $k-1$ dimensional projective subspace $l \subset \mathbb{C}
\mathbb{P}^n$ is moved by the action to another $k-1$ dimensional projective subspace. Therefore one gets the induced $T^n$ - action on the Grassmannian ${\rm Gr}(k, n+1)$.
For the projective space $\mathbb{P} (\wedge^k \mathbb{C}^{n+1})$ one has the induced K\"{a}hler structure, the induced Pl\"{u}cker coordinates $w_{i_1, ..., i_k}$, and therefore
the induced moment maps $F(\mu_i)$ must have  the explicit expressions. Indeed, they  read as follows
$$
F(\mu_i) = \frac{\sum_{(i_1, ..., i_k)} \delta (i, i_1, ..., i_k) \vert w_{i_1, ..., i_k} \vert^2}{\sum_{(i_1, ..., i_k)} \vert w_{i_1, ..., i_k} \vert^2},
$$
where symbol $\delta (i, i_1, ..., i_k)$ equals 1 if $i = i_j$ for certain $j$ or to zero otherwise.  

The geometrical meaning of the induced moment maps $F(\mu_i)$ is rather simple: the value of $F(\mu_i)$ at subspace $L \subset \mathbb{C}^{n+1}$ equals to
 the norm of the orthogonal projection of basis unit vector $v_i$ to $L$;
on the projective level the value is a derivation from the distance between point $[0: ... : 1: ...0]$ where 1 is on the $i$th place and the corresponding
projective subspace $l = \mathbb{P}(L)$. Essentially $F(\mu_i)$ measures the angle of $L$ to $v_i$.
 
 From this description it follows that $F(\mu_i)$ has two critical values 0 and 1, and for any other value $0< c_i <1$
the Hamiltonian vector field $X_{F(\mu_i)}$ does not vanish on the level set $N(c_i) = \{F(\mu_i) = c_i \} \subset {\rm Gr}(k, n+1)$. The critical values
0 and 1 corresponds to the following "ends": recall that ${\rm Gr}(k, n+1)$ under the choice of a vector $v \in \mathbb{C}^{n+1}$ can be decomposed into two parts
$$
{\rm Gr} (k, n+1) = {\rm Gr}(k-1, n+1) \cup {\rm tot} (E \to {\rm Gr}(k, n))
\eqno (1)
$$
 where $E$ is a vector bundle of rank $k$. The first part corresponds to subspaces $L$ which contain
$v$, --- and in our case it is equivalent to the fact that the norm of the  projection of $v$ to $L$ equals to 1; the zero section of $E$ corresponds to the subspaces $L$
such that the projection equals to 0, hence the zero section of $E$ consists of the subspaces contained by the orthogonal complement $< v>^{\perp}$. Therefore
the critical subsets of $F(\mu_i)$ can be described as follows: take decomposition (1) for vector $v= v_i$, then the first part of (1) forms the critical subset
with critical value 1, and the zero set ${\rm Gr}(k, <v_i>^{\perp})$ corresponds to critical subset with critical value 0.

On the other hand we know that the real part ${\rm Gr}_{\mathbb{R}}(k, n+1) \subset {\rm Gr}(k, n+1)$ exists, has right dimension and is transversal
to the toric action of each $F(\mu_i)$. Indeed, for any real $l \subset \mathbb{R} \mathbb{P}^n \subset \mathbb{C} \mathbb{P}^n$
the flow $\phi^t_{X_{\mu_i}}$ is either trivially acts on $l$ or it moves $l$ outside of $\mathbb{R} \mathbb{P}^n$. The trivial action
corresponds to the cases $F(\mu_i)(l) = 0$ or 1; otherwise $l$ must contain a real point with non trivial $i$'th - coordinate, the flow
scales this coordinate by $e^{it}$ and does not change the resting coordinates ---  hence the point must leave $\mathbb{R} \mathbb{P}^n$ (but
it comes there again for a moment  when $t = \pi$ which we will exploit below).

{\bf Mironov submanifolds of homogeneity 1.} As we have seen above complex Grassmanian ${\rm Gr} (k, n+1)$ is an algebraic variety which possesses the properties
one needs to apply the construction of Mironov cycles. At the same time instead of full $T^n$ - action, spanned by all the moment maps $F(\mu_i)$, one can reduce
the story to any subtorus $T^k$. To distinguish the cases we say that a Mironov cycle $T^k(S_{\mathbb{R}}(c_1, ..., c_k))$ has homogeneity $k$
if it is constructed using $k$ moment maps, derived from the complete set $(F(\mu_1), ..., F(\mu_n))$. Below we present an example
of Mironov submanifold constructed using a single moment map, say, $F(\mu_n)$.

Fix a non critical value $c_n \in (0; 1)$ and study first of all the restricted level set $S_{\mathbb{R}}(c_n) = N(c_n) \cap {\rm Gr}_{\mathbb{R}} (k, n+1)$.   
Since $c_n \neq 0$ we can exclude
the first component in the decomposition (1). Recall the description of the second part. The bundle $E_{\mathbb{R}} \to {\rm Gr}_{\mathbb{R}}(k, n)$ has as the base
(and the zero section) the space of $k$ - dimensional subspaces in the orthogonal complement $< v_n>^{\perp} \subset \mathbb{R}^{n+1}$. For any such subspace $L_0
\subset < v_n>^{\perp}$ the fiber $E_{\mathbb{R}}|_{L_0}$ consists of $k$ dimensional subspaces of the direct sum $\mathbb{R} <v_n> \oplus L_0$ which do not
contain $v_n$ (therefore they form $k$ - dimensional vector space).

Now, since we study $N(c_n)$, consider $k$ - dimensional subspace $L$ of $\mathbb{R} <v_n> \oplus L_0$ such that the projection of $v_n$ to $L$ is fixed;
it is equivalent to the condition that the angle between $v_n$ and $L$ is fixed. Since $c_n \neq 0$ such $L$ never coincides with $L_0$, and the intersection
$L \cap L_0 = M \subset L_0$ is a proper $k-1$ - dimensional subspace of $L_0$. Then it is not hard to see that if $M \subset L_0$ is fixed there exist exactly
two choices of such $k$ - dimensional subspaces in $\mathbb{R} <v_n> \oplus L_0$ with fixed angle, which contains $M$.

These arguments imply that the restricted level set $S_{\mathbb{R}}(c_n) \subset {\rm Gr}_{\mathbb{R}}(k, n+1)$ is isomorphic to the following manifold.
Take ${\rm Gr}_{\mathbb{R}}(k, n)$, take the tautological bundle $\tau \to {\rm Gr}_{\mathbb{R}}(k, n)$, take the dual bundle $\tau^*$
and at last take the "spherization" $S^{k-1}(\tau^*)$ of this bundle: then the total space
$$
{\rm tot} ( S^{k-1}(\tau^*) \to {\rm Gr}_{\mathbb{R}}(k, n))
\eqno (2)
$$
is isomorphic to $S_{\mathbb{R}}(c_n) = N(c_n) \cap {\rm Gr}_{\mathbb{R}}(k, n+1)$.

 Indeed, as we have seen above, over a point $[L_0] \in {\rm Gr}_{\mathbb{R}}(k, n)$ two $k$ - dimensional subspaces $L_1$ and $L_2$
are uniquely defined by the fixed angle to $v_n$ and the intersection $M = L_0 \cap L_1 = L_0 \cap L_2$. The intersection $M$ is given by
a point of the projectivization of the fiber $\tau^*|_{[L_0]}$, therefore one has the double covering of $\mathbb{P}(\tau^*)$ which is
the spherization of the fiber; globalization of the local picture leads to the answer, given in (2).

Further, following the strategy, we switch  on the Hamiltonian action generated by $X_{F(\mu_n)}$. Geometrically this means that the unit vector $v_n$ varies
in the family $\{e^{it} v_n \}$, and the fiber subspaces $L_1$ and $L_2$ vary as well in the space $\mathbb{C}^{n+1}$. However under the process $L_0$ stays stable
since it is contained by $<v_n>^{\perp}$ which is stable for the Hamiltonian action being the critical subset. At the same time it is easy to see that 
the rotation $\phi_{X_{F(\mu_n)}}^{t}$ interchanges $L_1$ and $L_2$ when $t = \pi$. 

The result of the $S^1$ -  action, generated by the moment map $F(\mu_n)$, on $S_{\mathbb{R}}(c_n)$ has been described by P. Pushkar' in paper [4]. 
Take the direct sum $\underline{S}^1 \times S^{k-1}(\tau^*) \to  {\rm Gr}_{\mathbb{R}}(k, n))$ of the trivial $S^1$ bundle and the spherization
of $\tau^*$; there one has the fiberwise diagonal action of $\mathbb{Z}_2$ given by simultaneous action of the standard antipodal involutions
on both the summands (note that both the summands are spheres). Factorizing with respect to this $\mathbb{Z}_2$ - action one gets the answer:
$$
S^1(S_{\mathbb{R}}(c_n)) = {\rm tot} (S^1 \times S^{k-1}(\tau^*))/\mathbb{Z}_2 \to  {\rm Gr}_{\mathbb{R}}(k, n)),
\eqno (3)
$$
where the fiber is the Pushkar submanifold $L_k \subset \mathbb{C}^k = \mathbb{C} \otimes \tau^*|_{[L_0]}$.
According to [4] (Proposition 1), we can characterize the topological type of the constructed Mironov cycle as follows:
 it is presented as a fiber bundle over
the real Grassmannian where the fiber is either $S^1 \times S^{k-1}$ for even $k$ or topologically non trivial $U(1)$ - bundle over $\mathbb{R} \mathbb{P}^{k-1}$ for odd  $k$
(this type was called {\it generalized Klein bottle} in [1]);
in both the cases the fiber bundle is topologically non trivial being associated with $\tau^* \to {\rm Gr}_{\mathbb{R}}(k, n)$.

In particular for  $k = 2$ the corresponding Mironov cycle is a $T^2$ - bundle: for ${\rm Gr}(2,3) = \mathbb{C} \mathbb{P}^2$ the construction gives
the standard Clifford torus since the base ${\rm Gr}_{\mathbb{R}}(2,2)$ is just a point; for ${\rm Gr} (2, 4)$ it gives a two - torus bundle over $\mathbb{R} \mathbb{P}^2$.

Note however  that $F(\mu_n)$ is rather simple moment map for the present case: one can study any integer valued linear combination of $F(\mu_i)$
as the moment map used for the construction of a homogeneity 1 Mironov cycle, and as we know from [2] the result can have different topological type.

We hope to continue the work in the future.

$$$$

{\bf References:}

[1] A. Mironov, {\it “New examples of Hamilton-minimal and minimal Lagrangian manifolds in $\mathbb{C}^n$ and $\mathbb{C} \mathbb{P}^n$”}, Sb. Math., 195:1 (2004) pp. 85–96;

[2] N. Tyurin, {\it "Lagrangian cycles of Mironov in algebraic varieties"}, submitted to Sb. Math;

[3] P. Griffits, J. Harris, {\it "Principles of algebraic geometry"}, NY, Wiley, 1978;

[4] P. Pushkar', {\it "Lagrange intersections in a symplectic space"},  Funct. An. and Its Appl., 34 (2000), pp. 288 - 292.

\end{document}